\documentclass[11pt]{article}
\begin{document}

\author{Oleg Pikhurko\thanks{Supported by a Senior Rouse Ball
 Studentship, Trinity College, Cambridge, UK.}\\
 DPMMS, Centre for Mathematical Sciences\\
 Cambridge University, Cambridge~CB3~0WB, England\\
 E-mail: {\tt O.Pikhurko@dpmms.cam.ac.uk}}

\newcommand{\eqref}[1]{\mbox{\rm(\ref{#1})}}
\newcommand{\pl}{}
\newcommand{\C}[1]{{\protect\cal #1}}
\newcommand{\binom}[2]{{#1\choose #2}}
\newcommand{\qed}{\nolinebreak\mbox{\hspace{5 true pt}%
\rule[-0.85 true pt]{3.9 true pt}{8.1 true pt}}}

\newtheorem{theorem}{Theorem}
\newtheorem{lemma}[theorem]{Lemma}
\newtheorem{corollary}[theorem]{Corollary}
\newtheorem{problem}[theorem]{Problem}

\newcommand{\Aut}{\mbox{\rm Aut}}
\newcommand{\cing}{\cong'}

\title{Foata's Bijection for Tree-Like Structures}
\maketitle

\begin{abstract} We present bijections enumerating $(k,m)$-trees,
$k$-gon trees, edge labelled $(2,1)$-trees, and other tree-like
structures. Our constructions are based on Foata's~\cite{foata:71}
bijection for cycle-free functions, which is simplified
here.\end{abstract}

\section{Introduction}

Here we enumerate certain tree-like structures by finding an explicit
bijection between the set in question and some trivially simple
set. Although some of the obtained formulas are known, the constructed
bijections are new results. Bijective proofs are of interest as they
carry more information than formulas alone: bijections allow us to
generate all objects one by one and it is often possible to enumerate
some subclasses by analysing the corresponding codes.

For example, the analytic formula for the number of vertex-labelled
$(m+1,m)$-trees was independently found by Beineke and
Pippert~\cite{beineke+pippert:69} and Moon~\cite{moon:69}. (We will
provide all definitions in due course.) Nevertheless, a number of
different bijective proofs of the result appeared as well,
see~\cite{renyi+renyi:70,foata:71,greene+iba:75,egecioglu+shen:88,chen:93}.

One of these papers, by Foata~\cite{foata:71}, contains an interesting
bijection for cycle-free functions. In Section~\ref{bijective:foata}
we describe Foata's construction with a simpler argument than the
original one.

Cycle-free functions provide a natural framework for encoding
different tree-like structures, in particular $(m+1,m)$-trees, as
demonstrated by Foata~\cite{foata:71}. Extending Foata's argument, we
present a bijection for general $(k,m)$-trees. (The analytic formula
was found by the author~\cite{pikhurko:99:jcta}.)

This method can enumerate some other objects. For example, we can find
a bijection for vertex labelled {\em $k$-gon trees} (also known as
{\em cacti} or {\em trees of polygons}), a structure that appears
in~\cite{chao+li:85,whitehead:88,peng:93,koh+teo:96}. The derived
formula for the number of vertex-labelled $k$-gon trees seems to be a
new result.

In order not to repeat the same portions of proof twice, we enumerate
a more general structure, {\em $H$-built-trees}, which includes
$(k,m)$-trees and $k$-gon trees as partial cases. Please refer to
Section~\ref{bijective:tree} for details.

In Section~\ref{tree:edge} we present a bijection for edge labelled
$(2,1)$-trees. This answers a question posed by
Cameron~\cite{cameron:95} which was motivated by the possibility that
such a bijection can simplify some of his enumeration results (or
proofs). However, although we answered Cameron's question, we were not
able to improve on~\cite{cameron:95}.

Unfortunately, we do not know any bijection (nor even a compact
analytic formula) giving the number of edge labelled $(k,m)$-trees for
$k\ge 3$.

\section{Foata's Bijection}\label{bijective:foata}

Given disjoint finite sets $A$, $B$, $C$ and a surjection $\gamma:B\to
A$, a function $f:A\to B\cup C$ is called {\em cycle-free} if for
every $b\in B$ the sequence $(f\circ\gamma)^i(b)$ eventually
terminates at some $c\in C$.

Foata~\cite[Theorem~1]{foata:71} exhibited a bijection between
$F(A,B,C,\gamma)$, the set of cycle-free functions, and the set of
functions $g:A\to B\cup C$ such that $g(a_1)\in C$ for some beforehand
fixed $a_1\in A$; this implies\begin{equation}\label{eq:\pl:foata}
|F(A,B,C,\gamma)|=|C|\big(|B|+|C|\big)^{|A|-1}.\end{equation}

We briefly describe a simpler construction than that
in~\cite{foata:71}. Fix some ordering of $A$. Let $f\in
F(A,B,C,\gamma)$. Let $Z=(z_1,\dots,z_s)$ denote the increasing
sequence of the elements in $A\setminus \gamma(f(A))$. (For convenience we
assume that $\gamma(c)=c$ for $c\in C$.) We build, one by one, $s$
sequences $\delta_1,\dots,\delta_s$ composed of elements in $B\cup C$.
Having constructed sequences $\delta_1,\dots,\delta_{i-1}$, let
$m_i\ge0$ be the smallest integer such that
$(\gamma\circ f)^{m_i+1}(z_i)$ either belongs to $C$ or occurs in at
least one of $\gamma\otimes\delta_1,\dots,\gamma\otimes
\delta_{i-1}$, where $\gamma\otimes
(x_1,\dots,x_i)\equiv(\gamma(x_1),\dots,\gamma(x_i))$. We define (mind the
order)\begin{equation}\label{eq:\pl:delta}
 \delta_i=\left( (f\circ\gamma)^{m_i}(f(z_i)),\
(f\circ\gamma)^{m_i-1}(f(z_i)),\ \dots,\ f(z_i)\right).\end{equation}

One can easily check that $Z$ is non-empty if $A$ is, every $m_i$
exists, and $\delta$, the juxtaposition product of the $s$ sequences
$\delta_1,\dots,\delta_s$, contains $|A|$ elements. (In fact, $\delta$
is but a permutation of $(f(a))_{a\in A}$.) The obtained sequence
$\delta$ of $|A|$ elements of $B\cup C$, which starts with an element
of $C$, corresponds naturally to the required function $g:A\to B\cup
C$.

Conversely, given $g$ (or $\delta$) we can reconstruct $Z=A\setminus
\gamma(g(A))$. Then, exactly $s=|Z|$ times, an element of
$\gamma\otimes \delta$ either belongs to $C$ or equals some preceding
element. These $s$ positions mark the beginnings of
$\delta_1,\dots,\delta_s$. Now we can restore the required $f$
by~\eqref{eq:\pl:delta}. To establish~\eqref{eq:\pl:foata} completely, one has to check
easy details.

\section{$H$-Built-Trees}\label{bijective:tree}

The following notion of {\em\ $(k,m)$-tree} was suggested independently
by Dewdney~\cite{dewdney:74} and Beineke and
Pippert~\cite{beineke+pippert:77}.

Let us agree that the vertex set is $[n]\equiv\{1,\dots,n\}$. Fix the
{\em edge size} $k$ and the {\em\ overlap size}
$m\in[0,k-1]\equiv\{0,\dots,k-1\}$. We refer to $k$-subsets and
$m$-subsets of $[n]$ as {\em edges} and {\em laps}\index{lap}
respectively. A non-empty $k$-graph (i.e.\ a $k$-uniform set system)
without isolated vertices is called a {\em\ $(k,m)$-tree} if we can
order its edges, say $E_1,\dots,E_e$, so that, for every $i\in[2,e]$,
there is $j\in[i-1]$ such that $|E_i\cap E_{j}|=m$ and $(E_i\setminus
E_{j})\cap \left(\cup_{h=1}^{i-1} E_h\right)=\emptyset$. In other
words, we start with a single edge and can consecutively affix a new
edge along an $m$-subset of an existing edge. For example, a
$(k,0)$-tree consists of disjoint edges; $(2,1)$-trees are usual
(Cayley) trees.

Adopting the ideas of Foata~\cite{foata:71}, we present a bijection
for $(k,m)$-trees. In fact, we enumerate a more general structure
defined as follows.

Let $H$ be any $m$-graph on $[k]$. An {\em\ $H$-built-tree}
$(T,\{H_1,\dots,H_e\})$ consists of a usual $(k,m)$-tree $T$ with edges
$E_1,\dots,E_e$ plus $H$-graphs $H_i$ on $E_i$, $i\in[e]$, such that
if $E_i\cap E_j$ is a lap (that is, has size $m$), then it is an edge
of both $H_i$ and $H_j$ for any distinct $i,j\in[e]$. Let
$n=e(k-m)+m$ denote the total number of vertices and let
 $$f=\left|\cup_{i\in[e]}E(H_i)\right|=e(l-1)+1,$$
 where $l=e(H)=|E(H)|$ denotes the number of edges in $H$. Also, let $\C 
R_H$ be the set of distinct $H$-graphs on $[k]$ rooted at $[m]$, that
is, containing $[m]$ as an edge. Clearly,
 $$|\C R_H|=\frac{k!l}{ |\Aut(H)|\binom km}.$$
 An $H$-built-tree is {\em rooted}\index{rooted tree} on
an $m$-set $L$ if $L\in \cup_{i\in[e]}E(H_i)$.

\begin{theorem}\label{th:\pl:bijective:general} There is a bijection between the set $Y$ of
$H$-built-trees on $[n]$ rooted at $[m]$ and the set
 $$ Z=F(A,B,C,\gamma)\times\prod_{i=1}^e \left(X_i \times \C R_H\right),$$
 where  $A=[e]$, $B=[e]\times[l-1]$, $C=\{[m]\}$, $\gamma$
is the coordinate projection $B\to A$, and
$X_{i}=\big[{(k-m)(e-i+1)-1\choose k-m-1}\big]$. In
particular,\begin{eqnarray*}
 |Y| &=& f^{e-1}\left(\frac{k!\,l}{\binom km |\Aut(H)|}\right)^e
\prod_{i=1}^e{(k-m)(e-i+1)-1\choose k-m-1}\\
 &=&\frac{(e(k-m))!\,f^{e-1}}{e!}
\left(\frac{m!\, l}{|\Aut(H)|}\right)^e.\end{eqnarray*}\end{theorem}
 \smallskip{\it Proof.}  Given an $H$-built-tree $T$ rooted at $[m]$, order its edges
$E_1,\dots,E_e$ so that $[m]\in E(H_1)$ and each $E_i$, $i\in[2,e]$,
shares a lap with some $E_j$, $j<i$.

Correspond an edge $E_i$ to the lap $g'(E_i)=E_i\cap \cup_{j=1}^{i-1}
E_j$, $i\in[2,e]$. (We agree that $g'(E_1)=[m]$.)  Call the set
$f(E_i)=E_i\setminus g'(E_i)$ the {\em free part} of $E_i$; the free
parts partition $[n]\setminus[m]$. Clearly, these definitions of $g'$
and $f$ do not depend on the particular ordering.

Relabel the edges by $D_1,\dots,D_e$ so that $d_i=\min f(D_i)$
increases; let $H_i'$ denote the corresponding $H$-graph on
$D_i$. Label, in the colex order, all edges (laps) of $H_i'$ but
$g'(D_i)\in E(H_i')$ by $(i,j)$, $j=1,\dots,l-1$. Note that now we
have indexing of the edges of $T$ by $A$, namely $(D_i)_{i\in A}$, and
of the laps of $T$ by $B\cup C$. Let $g:A\to B\cup C$ be the map
corresponding to $g'$. A moment's thought reveals that $g$ is
cycle-free.

Repeat the following for $i=1,\dots,e$. Index, in the colex order, the
$(k-m-1)$-subsets of $(\cup_{j=i}^e f(D_j))\setminus\{d_i\}$ by the
elements of $X_i$ and let $x_i\in X_i$ be the index corresponding to
$f(D_i)\setminus\{d_i\}$. Consider the bijection $h:D_i\to [k]$ such
that $h$ is monotone on $g'(D_i)$ and $f(D_i)$ which are respectively
mapped onto $[m]$ and $[m+1,k]$. Let $R_i\in\C R_H$ be the
image of $H_i'$ under $h$.

Now, $(g,x_1,R_1,\dots,x_e,R_e)\in Z$ is the `code' of $T\in Y$.

Conversely, given an element $(g,x_1,R_1,\dots,x_e,R_e)\in Z$ we can
consecutively reconstruct the sequence $(d_i,f(D_i))$, $i=1,\dots,e$.
Indeed, $d_i$ is the smallest element of $V=[n]\setminus
\big((\cup_{j=1}^{i-1} f(D_j))\cup[m]\big)$ while
$f(D_i)\setminus\{d_i\}$ is the $x_i$th $(k-m-1)$-subset of
$V\setminus\{d_i\}$. For $i\in A$ with $g(i)\in C$, we have
$D_i=[m]\cup f(D_i)$ and (knowing $g'(D_i)=[m]$ and $f(D_i)$), we can
determine $H_i'$ from $R_i$; then we can recover the lap corresponding
to $(i,j)\in B$ as the $j$th element of $E(H_i')\setminus\{[m]\}$,
$j\in[l-1]$.

Likewise, we can reconstruct all information about $D_i$
for any $i\in A$ with $g(i)$ being already associated with a lap. As
$f$ is cycle-free, all edges are eventually identified, producing $T\in
Y$.

The plain verification shows that we have indeed a bijective
correspondence between $Y$ and $Z$.\qed \medskip

It is easy to see that $(k,m)$-trees bijectively correspond to
$K_k^m$-built-trees on the same vertex set, where $K_k^m$ denotes the
complete $m$-graph of order $k$. Now, $n=e(k-m)+m$, $l=\binom km$,
$f=e(l-1)+1$, $|\C R_{K_K^m}|=1$, and we deduce the
following.

\begin{corollary}\label{cr:\pl:kmtree} The number of $(k,m)$-trees on $[n]$ rooted at
$[m]$ equals\begin{equation}\label{eq:\pl:1}
 R_{km}(e)=\frac{(n-m)!f^{e-1}}{e!\left((k-m)!\right)^e}.\qed\end{equation}\end{corollary}

\smallskip\noindent{\bf Remark.}  Clearly, the number of vertex labelled $(k,m)$-trees with $n$
vertices is equal to ${n\choose m} R_{km}(e)/f$, which gives
precisely~\cite[formula~(1)]{pikhurko:99:jcta}.\smallskip

As another consequence, we can enumerate {\em\ $k$-gon trees}, which are
inductively defined as follows. A $k$-gon (that is, $C_k$, a cycle of
length $k$) is a $k$-gon tree. A $k$-gon tree with $e+1$ $k$-gons is
obtained from a $k$-gon tree with $e$ $k$-gons by adding $k-2$ new
vertices and a new $k$-gon through these vertices and an already
existing edge.

It is trivial to check that if a union of $C_k$-graphs can be formed
into a $C_k$-built-tree, then the latter is uniquely defined.
Thus, $k$-gon trees are in bijective correspondence with
$C_k$-built-trees. We have $m=2$, $n=e(k-2)+2$, $f=e(k-1)+1$, and $|\C
R_{C_k}|=(k-2)!$, so we obtain that there are $(e(k-2))!\, f^{e-1}/e!$
rooted $k$-gon trees with $e$ $k$-gons, which implies the following
result.

\begin{corollary}\label{cr:\pl:k-gon} The
number of vertex labelled $k$-gon trees with $e$ $k$-gons is
 $$\frac{(e(k-2)+2)!\,(e(k-1)+1)^{e-2}}{2 (e!)},\quad k\ge 3.\qed$$\end{corollary}

\section{Edge Labelled Trees}\label{tree:edge}

Cameron~\cite{cameron:95} enumerates certain classes of what is called
there {\em two-graphs}: reduced, $5$-free, and $(5,6)$-free. Please
refer to his work for all definitions and details. Also, he defines,
for a given (Cayley) tree $T$, the equivalence relation $\cong$ on its
edges which is the smallest one such that two edges are related if
they intersect at a vertex of degree $2$ in $T$. For example, $T$ is
series-reduced if and only if $\cong$ is the identity relation.

Cameron had to count the number $S_n$ of trees
with $n$ labelled edges when we do not distinguish trees obtained by
permuting labels within the $\cong$-classes. He found the following
formula~\cite[Proposition~3.5(a)]{cameron:95}:\begin{equation}\label{eq:\pl:cameron}
 S_n=\sum_{k=1}^n S(n,k)\frac{1}{k+1} \sum_{j=0}^{k-1}
(-1)^j\binom{k+1}j \binom{k-1}j  j!(k-j+1)^{k-j-1},\end{equation}
 where $S(n,k)$ is the Stirling number of the second kind. The sequence
$(S_n)$ starts as $1,1,2,8,52,\dots$ and probably cannot be
represented in a closed expression but, of course, one can try to
simplify~\eqref{eq:\pl:cameron}. Cameron~\cite{cameron:95} asks the following question.

\begin{problem}[Cameron]\label{pr:\pl:cameron} Describe a constructive
bijection for edge labelled trees, not
going via vertex labelling. Describe the equivalence relation $\cong$
in terms of this code.\end{problem}

The motivation for the problem was that such a code might
simplify~\eqref{eq:\pl:cameron}. Although we answer here this question, we do
not see how our bijection can simplify~\eqref{eq:\pl:cameron} or its proof
from~\cite{cameron:95}.

Let us describe our construction.  We use Foata's~\cite{foata:71}
bijection for cycle-free functions.

Let $e_1,\dots,e_n$ be the edges. Suppose $e_1=\{a,b\}$; this edge
will play a special role. Let $A=B=\{e_2,\dots,e_{n}\}$, $C=\{a,b\}$
and $\gamma:B\to A$ be the identity function. Let us correspond an $f\in
F(A,B,C,\gamma)$ to a given tree $T$.  Each
edge $e$ can be connected to $e_1$ by the unique path in $T$. If $e$ is
incident to $e_1$, then let $f(e)$ be equal to their common vertex;
otherwise, let $f(e)$ be the first edge on the path from $e$ to
$e_1$. This gives a correspondence between twice the number of
edge-labelled trees (we can label the two vertices of $e_1$ by $a$ and
$b$ in two different ways) and $F(A,B,C,\gamma)$.  Foata's bijection
shows that $|F(A,B,C,\gamma)|=|C|(|A|+|C|)^{|A|-1}$, which implies, as
desired, that the number of edge-labelled trees with $n$ edges is
\mbox{$(n+1)^{n-2}$}.

To make this correspondence one-to-one, we can consider only a half of
the codes, e.g.\ those starting with $a$. They correspond to trees in
which the path from the leaf $e_i$ with the smallest index $i\in[2,n]$
to $e_1$ hits $e_1$ at $a$.

Of course, the code is rather simple; we describe briefly only
one direction. A code $\delta$ is a sequence of length $n-1$
consisting of elements in $\{a,b,e_2,\dots,e_{n}\}$ and staring with
$a$. The set $Z\subset\{e_2,\dots,e_n\}$ of edges which do not
occur in the sequence consists of leaves. (If $b$ does not
occur, then $e_1$ is also a leaf.)  Clearly, an element of $\delta$
equals either $a$ or $b$ or some previously occurring element exactly
$s=|Z|$ times. Cut $\delta$ before each such element; we have $s$
pieces $\delta_1,\dots,\delta_s$. Append the $i$th element $z_i$ of
$Z$ to the end of $\delta_i$ to obtain $\delta_i'$, $i\in[s]$.

The reversed sequence $\delta_i'$ describes the initial segment
$P_i'$ of the path $P_i$ from the element $z_i\in Z$ to $e_1$
until it hits $e_1$ or some previous path $P_j$, $i\in[s]$. Clearly,
this determines some tree.

How can we read the $\cong$-relation from $\delta$? First, let $\cing$
be the minimal equivalence relation on $\{a,b,e_2,\dots,e_n\}$ such
that $e_i\cing e_j$ if $e_i$ and $e_j$ intersect at a vertex of degree
$2$, $2\le i<j\le n$, and $x\cing e_i$ if $x$ is a degree-$2$ vertex
incident to $e_i$, $x\in\{a,b\}$, $i\in[2,n]$. (Informally, we cut
$e_1$ in the middle and take the usual $\cong$-relation on the both
created components separately.) Clearly, $\cong$ is obtained from
$\cing$ by identifying $a$ and $b$ into a single element $e_1$, so let
us indicate how to determine the latter relation.

Take any maximal contiguous subsequence $S\subset \delta$ consisting
of elements that occur in $\delta$ exactly once. Clearly, $S$ lies
entirely within some $\delta_i$ and $S\cup\{y\}$ is a
$\cing$-equivalence class, where $y$ is the symbol following $S$ in
$\delta_i'$. Conversely, it is easy to check that all non-trivial
$\cing$-classes are obtained this way, as required.


\begin{thebibliography}{10}

\bibitem{beineke+pippert:69}
L.~W. Beineke and R.~E. Pippert.
\newblock The number of labelled {$k$}-dimensional trees.
\newblock {\em J.\ Combin.\ Theory}, 6:200--205, 1969.

\bibitem{beineke+pippert:77}
L.~W. Beineke and R.~E. Pippert.
\newblock On the structure of {$(m,n)$}-trees.
\newblock In {\em Proc.\ 8th {Southeast} {Conference} on {Combinatorics},
  {Graph Th.}\ and {Computing} {(Louisiana State Univ., Baton Rouge)}}, volume
  XIX of {\em Congressus Numeratium}, pages 75--80. Utilitas Math., Winnipeg,
  1977.

\bibitem{cameron:95}
P.~J. Cameron.
\newblock Counting two-graphs related to trees.
\newblock {\em Electronic J.\ Combin.}, 2(\#R4):8pp, 1995.

\bibitem{chao+li:85}
C.-Y. Chao and N.-Z. Li.
\newblock Trees of polygons.
\newblock {\em Archiv Math.}, 45:180--185, 1985.

\bibitem{chen:93}
W.~Y.~C. Chen.
\newblock A coding algorithm for {R\'enyi} trees.
\newblock {\em J.\ Combin.\ Theory\ {\rm (A)}}, 63:11--25, 1993.

\bibitem{dewdney:74}
A.~K. Dewdney.
\newblock Multidimensional tree-like structures.
\newblock {\em J.\ Combin.\ Theory\ {\rm (B)}}, 17:160--169, 1974.

\bibitem{egecioglu+shen:88}
O.~Egecioglu and L.-P. Shen.
\newblock A bijective proof for the number of labelled {$q$}-trees.
\newblock {\em Ars Combinatoria}, 25:3--30, 1988.

\bibitem{foata:71}
D.~Foata.
\newblock Enumerating {$k$}-trees.
\newblock {\em Discrete Math.}, 1:181--186, 1971.

\bibitem{greene+iba:75}
C.~Greene and G.~A. Iba.
\newblock Caley's formula for multidimensional trees.
\newblock {\em Discrete Math.}, 14:1--11, 1975.

\bibitem{koh+teo:96}
K.~M. Koh and C.~P. Teo.
\newblock Chromaticity of series-parallel graphs.
\newblock {\em Discrete Math.}, 154:289--295, 1996.

\bibitem{moon:69}
J.~W. Moon.
\newblock The number of labeled {$k$}-trees.
\newblock {\em J.\ Combin.\ Theory}, 6:196--199, 1969.

\bibitem{peng:93}
Y.-H. Peng.
\newblock On the chromatic uniqueness of certain trees of polygons.
\newblock {\em J. Austral. Math. Soc. Ser. A}, 55:403--410, 1993.

\bibitem{pikhurko:99:jcta}
O.~Pikhurko.
\newblock Enumeration of labelled {$(k,m)$}-trees.
\newblock {\em J.\ Combin.\ Theory\ {\rm (A)}}, 86:197--199, 1999.

\bibitem{renyi+renyi:70}
C.~{R\'enyi} and A.~{R\'enyi}.
\newblock The {Pr\"ufer} code for {$k$}-trees.
\newblock In P.~{Erd\H os} et~al., editors, {\em Combinatorial Theory and Its
  Applications}, pages 945--971. North-Holland, Amsterdam, 1970.

\bibitem{whitehead:88}
E.~G. Whitehead, Jr.
\newblock Chromatic polynomials of generalized trees.
\newblock {\em Discrete Math.}, 72:391--393, 1988.

\end{thebibliography}
\end{document}